\documentclass[11pt]{amsart}

\numberwithin{equation}{section}

\usepackage{amsfonts}
\usepackage{amssymb}
\newtheorem{theorem}{Theorem}[section]
\newtheorem{corollary}[theorem]{Corollary}
\newtheorem{lemma}[theorem]{Lemma}
\newtheorem{proposition}[theorem]{Proposition}

\theoremstyle{definition}

\newtheorem{remark}[theorem]{Remark}

\newcommand{\Ker}{\operatorname{Ker}}
\newcommand{\dom}{\operatorname{Dom}}

\begin{document}

\title[Aspects of the $L^{2}$-Sobolev theory of the $\overline{\partial}$-Neumann problem]{Aspects of the $L^{2}$-Sobolev theory of the $\overline{\partial}$-Neumann problem}

\author{Emil J. Straube}
\address{Department of Mathematics\\
Texas A\&M University\\
College Station, Texas, 77843--3368}
\email{straube@math.tamu.edu}

\thanks{Research supported in part by NSF grant DMS 0500842 and by a Texas A\&M University Faculty Development Leave}


\begin{abstract}
The $\overline{\partial}$-Neumann problem is the fundamental boundary value problem in several complex variables. It features an elliptic operator coupled with non-coercive boundary conditions. The problem is globally regular on many, but not all, pseudoconvex domains.

We discuss several recent developments in the $L^{2}$-Sobolev theory of the $\overline{\partial}$-Neumann problem that concern compactness and global regularity.
\end{abstract}

\maketitle

\section{Introduction}

The $\overline{\partial}$-Neumann problem was formulated in the fifties by D.~C.~Spencer as a means to generalize the theory of harmonic integrals (i.e. Hodge theory) to non-compact complex manifolds. For domains in $\mathbb{C}^{n}$, which is the context we will restrict ourselves to almost exclusively in this paper, the problem can be formulated as follows. Denote by $\Omega$ a pseudoconvex domain in $\mathbb{C}^{n}$, and by $L_{(0,q)}^{2}(\Omega)$ the space of $(0,q)$-forms on $\Omega$ with square integrable coefficients. Each such form can be written uniquely as a sum 
\begin{equation}\label{intro1}
u = \sum^{\prime}_{J}u_{J}d\overline{z_{J}} \;, 
\end{equation}
where $J=(j_{1}, \cdots , j_{q})$ is a multi-index with $j_{1} < j_{2} < \; \cdots \; < j_{q}$, $d\overline{z_{J}} = d\overline{z_{j_{1}}}\wedge\; \cdots \; \wedge d\overline{z_{j_{q}}}$, and the $^{\prime}$ indicates summation over  increasing multi-indices. The inner product
\begin{equation}\label{intro2}
 (u, v) = \left (\sum^{\prime}_{J}u_{J}d\overline{z_{J}},\sum^{\prime}_{J}v_{J}d\overline{z_{J}} \right ) = \sum^{\prime}_{J}\int_{\Omega}u_{J}\overline{v_{J}}  
\end{equation}
turns $L^{2}_{(0,q)}(\Omega)$ into a Hibert space. Set
\begin{equation}\label{intro3}
 \overline{\partial}\left(\sum^{\prime}_{J}u_{J}d\overline{z_{J}}\right) = \sum_{j=1}^{n}\sum^{\prime}_{J}\frac{\partial u_{J}}{\partial \overline{z_{j}}}d\overline{z_{j}}\wedge d\overline{z_{J}} \;, 
\end{equation}
where the derivatives are computed as distributions, and the domain of $\overline{\partial}$ is defined to consist of those $u \in L^{2}_{(0,q)}(\Omega)$ where the result is a $(0,q+1)$-form with square integrable coefficients. Then $\overline{\partial} =  \overline{\partial}_{q}$ is a closed, densely defined operator from $L^{2}_{(0,q)}(\Omega)$ to $L^{2}_{(0,q+1)}(\Omega)$, and as such has a Hilbert space adjoint. This adjoint is denoted by $\overline{\partial}_{q}^{*}$. (We will not use the subscripts when the form level at which the operators act is clear or not an issue.) 
One can check that $\overline{\partial}\overline{\partial}=0$, so that we arrive at a complex, the $\overline{\partial}$ (or Dolbeault)-complex:
\[ L^{2}(\Omega) \stackrel{\overline{\partial}}{\rightarrow} L^{2}_{(0,1)}(\Omega)
\stackrel{\overline{\partial}}{\rightarrow} L^{2}_{(0,2)}(\Omega) 
\stackrel{\overline{\partial}}{\rightarrow}
\; \cdots \; \stackrel{\overline{\partial}}{\rightarrow} L^{2}_{(0,n)}(\Omega) \stackrel{\overline{\partial}}{\rightarrow} 0 \; \;.\]
In analogy to the Laplace-Beltrami operator associated to the DeRham complex on a Riemannian manifold, one forms the complex Laplacian
\begin{equation}\label{intro4}
 \Box_{q} = \overline{\partial}_{q-1}\overline{\partial}_{q-1}^{*} + \overline{\partial}_{q}^{*}\overline{\partial}_{q} \;, 
\end{equation}
with domain so that the compositions are defined. The $\overline{\partial}$-complex is elliptic.
The $\overline{\partial}$-Neumann problem is the problem of inverting $\Box_{q}$; that is, given $v \in L^{2}_{(0,q)}(\Omega)$, find $u \in \dom(\Box_{q})$ such that $\Box_{q}u = v$. Note that $\dom(\Box_{q})$ involves the two boundary conditions $u \in \dom(\overline{\partial}^{*})$ and $ \overline{\partial}u \in \dom(\overline{\partial}^{*})$; these are the $\overline{\partial}$-Neumann boundary conditions. The condition $u \in \dom(\overline{\partial}^{*})$ is equivalent to a Dirichlet condition for the (complex) normal component of $u$. Similarly, the condition $\overline{\partial}u \in \dom(\overline{\partial}^{*})$ is equivalent to a Dirichlet condition on the normal component of $\overline{\partial}u$, that is, a \emph{complex (or $\overline{\partial}$-) Neumann condition} for $u$.

From the point of view of partial differential equations, the $\overline{\partial}$-Neumann problem represents the prototype of a problem where the operator is elliptic, but the boundary conditions are not coercive, so that the classical elliptic theory does not apply. From the point of view of several complex variables, the importance of the problem stems from the fact that its solution provides a Hodge decomposition in the context of the $\overline{\partial}$-complex, together with the attendant elegant machinery (as envisioned by Spencer). For example, such a decomposition readily produces a solution to the inhomogeneous $\overline{\partial}$ equation, as follows. Assume for the moment that $\Box_{q}$ has a (bounded) inverse in $L^{2}_{(0,q)}(\Omega)$, say $N_{q}$. Then we have the orthogonal decomposition
\begin{equation}\label{intro6}
u = \overline{\partial}\overline{\partial}^{*}N_{q}u + \overline{\partial}^{*}\overline{\partial}N_{q}u \;, \; u \in L^{2}_{(0,q)}(\Omega) \; .
\end{equation}
If $\overline{\partial}u =0$, then $\overline{\partial}^{*}\overline{\partial}N_{q}u$ is $\overline{\partial}$-closed as well (from \eqref{intro6}). Consequently,  $\overline{\partial}^{*}\overline{\partial}N_{q}u = 0$ (since it is also orthogonal to $\Ker(\overline{\partial})$), and
\begin{equation}\label{intro7}
u = \overline{\partial}\left(\overline{\partial}^{*}N_{q}u\right ) \; , 
\end{equation}
with $\|\overline{\partial}^{*}N_{q}u\|^{2} = (\overline{\partial}\overline{\partial}^{*}N_{q}u, N_{q}u) \leq C\|u\|^{2}$. Thus the operator $\overline{\partial}^{*}N$ provides an $L^{2}$-bounded solution operator to $\overline{\partial}$. In fact, this operator gives the (unique) solution orthogonal to $\Ker(\overline{\partial})$ (equivalently: the solution with minimal norm). This solution is called the canonical (or Kohn) solution.

That $\Box_{q}$ \emph{does} have a bounded inverse $N_{q}$ on bounded pseudoconvex domains was known by the mid sixties. Kohn (\cite{Kohn61, Kohn63,Kohn64}) solved the $\overline{\partial}$-Neumann problem for strictly pseudoconvex domains, showing that in this case, not only is there an $L^{2}$-bounded inverse, but $N_{q}$ exhibits a subelliptic gain of one derivative as measured in the $L^{2}$-Sobolev scale. Another interesting approach was given by Morrey in \cite{Morrey64}. H\"{o}rmander (\cite{Hormander65}, see also Andreotti-Vesentini \cite{AV65} for similar techniques) proved certain Carleman type estimates which in the case of bounded pseudoconvex domains imply the existence of $N_{q}$ as a bounded self-adjoint operator on $L^{2}_{(0,q)}(\Omega)$. Early applications included embedding of real analytic manifolds (\cite{Morrey58, Morrey64}), a new solution of the Levi problem (\cite{Kohn63}), a new proof of the Newlander-Nirenberg theorem on integrable almost complex manifolds (\cite{Kohn63}), and in general, an approach to several complex variables which takes advantage of the then newly developed $\overline{\partial}$-methods (\cite{Hormander65, Hormander73}). Interesting `eyewitness' accounts of this foundational period by two of the principals appear in \cite{Hormander03} and \cite{Kohn04}, respectively.

It is not hard to see that Kohn's results for strictly pseudoconvex domains are optimal: $N$ can never gain more than one derivative, and it can gain one derivative only when the domain is strictly pseudoconvex. However, under what circumstances subellipticity with a fractional gain of less than one derivative holds was not understood until the early eighties. Kohn gave sufficient conditions, satisfied for example when the boundary is real-analytic, in \cite{Kohn79}. In deep work, his students Catlin (\cite{Catlin83, Catlin84a, Catlin87b}) and D'Angelo (\cite{DA79, DA80, DA82}) resolved this question: on a smooth bounded pseudoconvex domain in $\mathbb{C}^{n}$, the $\overline{\partial}$-Neumann problem is subelliptic if and only if each boundary point is of finite type,that is, the order of contact, at the point, of complex varieties with the boundary is finite. For more on these ideas, see \cite{D'Angelo93, D'AngeloKohn99, D'Angelo03, Kohn04b}.

When $N_{q}$ does not gain derivatives, but is still compact (as an operator on $L^{2}_{(0,q)}(\Omega)$), it follows from work of Kohn and Nirenberg (\cite{KohnNirenberg65}) that $N_{q}$ preserves the Sobolev spaces $W^{s}_{(0,q)}(\Omega)$ for all $s \geq 0$. In particular, $N_{q}$ preserves $C^{\infty}_{(0,q)}(\overline{\Omega})$ (it is globally regular). Work of Catlin  (\cite{Catlin84b}, compare also Takegoshi \cite{Takegoshi91}) and Sibony (\cite{Sibony87b}) shows that compactness provides indeed a viable route to global regularity: the compactness condition can be verified on large classes of domains. We refer the reader to \cite{FuStraube99} for a survey of compactness in the $\overline{\partial}$-Neumann problem. In section 2 below, we will discuss some developments that have occurred since the publication of \cite{FuStraube99}.

In addition to the `usual' (pde) reasons for studying regularity properties of a differential operator, there are, in the case of the $\overline{\partial}$-Neumann problem, the implications global regularity of the $\overline{\partial}$-Neumann operator has for several complex variables. Chief among these is the relevance for boundary behavior of biholomorphic or proper holomorphic maps. Namely, if $\Omega_{1}$ and $\Omega_{2}$ are two bounded pseudoconvex domains in $\mathbb{C}^{n}$ with smooth boundaries, such that the $\overline{\partial}$-Neumann operator on $(0,1)$-forms (i.e. $N_{1}$) on $\Omega_{1}$ is globally regular, then any proper holomorphic map from $\Omega_{1}$ to $\Omega_{2}$ extends smoothly to the boundary of $\Omega_{1}$ (\cite{BC82, DiederichFornaess82}). This is a highly nontrivial result: in contrast to the one variable situation, where the result is classical, the general case in higher dimensions, even for biholomorphic maps, is open. An exposition of the ideas and issues involved here can be found in \cite{Bell90, Krantz92, ChenShaw01}.

That global regularity holds on large classes of pseudoconvex domains where local regularity or even compactness fail was shown in the early nineties by Boas and the author (\cite{BoasStraube91, BoasStraube93}). They proved in \cite{BoasStraube91} that if $\Omega$ admits a defining function whose complex Hessian is positive semi-definite at points of the boundary (a condition slightly more restrictive than pseudoconvexity), then the $\overline{\partial}$-Neumann problem is globally regular (for all $q$). This class of domains includes in particular all (smooth) convex domains. The proof is based on the existence of certain families of vector fields which have good approximate commutator properties with $\overline{\partial}$. In \cite{BoasStraube93}, the authors studied the situation when the boundary points of infinite type form a complex submanifold (with boundary) of the boundary of the domain. They identified a DeRham cohomology class associated to the submanifold as the obstruction to the existence of the vector fields needed. In particular, a simply connected complex manifold in the boundary is benign for global regularity of the $\overline{\partial}$-Neumann problem. It is noteworthy that this cohomology class also plays a role in deciding whether or not the closure of the domain admits a Stein neighborhood basis (\cite{BF78}).

The question whether global regularity holds on general pseudoconvex domains turned out to be very difficult and was resolved only in the mid nineties. Barrett (\cite{Barrett92}, see also \cite{Barrett84} and \cite{Kisel91} for predecessors) showed that on the worm domains of Diederich and Forn\ae ss (\cite{DF77}), $N_{1}$ does not preserve $W^{s}_{(0,1)}(\Omega)$ for $s$ sufficiently large, depending on the winding (that is, exact regularity fails). Christ (\cite{Christ96}, see also \cite{ Christ98, Christ99}) resolved the question by proving certain a priori estimates for $N_{1}$ on these domains that would imply exact regularity in Sobolev spaces (and thus would contradict Barrett's result) if $N_{1}$ were to preserve the space of forms smooth up to the boundary.  

In section 3, we will discuss some recent developments. In \cite{StraubeSucheston03} and \cite{ForstnericLaurent05}, the authors consider the case where the boundary is finite type except for a Levi-flat piece which is `nicely' foliated by complex hypersurfaces. Whether or not the families of vector fields with good approximate commutator properties with $\overline{\partial}$ exist turns out to be equivalent to a property of the Levi foliation much studied in foliation theory, namely whether or not the foliation can be defined globally by a closed one-form. Sucheston and the author showed in \cite{StraubeSucheston02} that the approaches via plurisubharmonic defining functions and vector fields with good approximate commutator properties with $\overline{\partial}$ are actually equivalent, when suitably reformulated. This left two main avenues to global regularity: compactness and plurisubharmonic defining functions and/or good vector fields. These two approaches were unified by the author in \cite{Straube05}, via a new sufficient condition for global regularity.

Detailed accounts of the $\overline{\partial}$-Neumann theory, from different points of view, may be found in \cite{FollandKohn72, Hormander73, Krantz92, ChenShaw01, LiebMichel02, Ohsawa02}. Developments up to about ten years ago are covered in \cite{BoasStraube99} (for compactness, see \cite{FuStraube99}). Two recent informative surveys that concentrate on topics not covered here are the following: \cite{Shaw05} deals with estimates on Lipschitz domains, and \cite{McNeal05} discusses applications obtained by inserting a so called twisting factor into the $\overline{\partial}$-complex.

Although this paper concentrates on aspects of the $\overline{\partial}$-Neumann problem on domains in $\mathbb{C}^{n}$, we wish to mention a very important and fruitful development, namely the application of $L^{2}$-methods to algebraic and complex geometry, and vice versa. For expositions of this very active area of research, we refer the reader to \cite{Demailly01, Demailly02, Siu02, Siu05}.

\section{Compactness}

It is of interest in several contexts to know whether or not the $\overline{\partial}$-Neumann operator is, or is not, compact. Examples include global regularity (\cite{KohnNirenberg65}), the Fredholm theory of Toeplitz operators(see e.g. \cite{HI97}), and existence or non-existence of Henkin-Ramirez type kernels for solving $\overline{\partial}$ (\cite{HeferLieb00}). 
A fairly comprehensive discussion of compactness, up to about 1999, is in \cite{FuStraube99}. There one also finds complete proofs and/or references for background material. In this section, we review some recent developments.

\vspace{0.1in}
\textbf{2.1 Sufficient conditions for compactness.} In \cite{Catlin84b} Catlin introduced a sufficient condition for compactness which he called property($P$). The boundary of a domain is said to satisfy property($P$) if for every positive number $M$ there are an open neighborhood $U_{M}$ of $b\Omega$ and a plurisubharmonic function $\lambda_{M} \in C^{2}(U_{M} \cap \Omega)$ with $0 \leq \lambda_{M} \leq 1$, such that for all $z \in U_{M} \cap \Omega$,
\begin{equation}\label{P}
\sum_{j,k}^{n}\frac{\partial^{2}\lambda_{M}}{\partial z_{j}\partial\overline{z_{k}}}(z)w_{j}\overline{w_{k}} \geq M|w|^{2} \; .
\end{equation}
(This is somewhat more general than the formulation given in \cite{Catlin84b}, but see also \cite{FuStraube99}.) Property($P$) can be very nicely reformulated, on Lipschitz domains, in the spirit of Oka's lemma: the boundary satisfies property($P$) if and only if it locally admits functions $\rho$ comparable to minus the boundary distance, such that the complex Hessian of $-\log(-\rho)$ tends to infinity upon approach to the boundary (\cite{Harrington05}). If the boundary of the bounded pseudoconvex domain satisfies property($P$), then the $\overline{\partial}$-Neumann operator $N_{q}$ on $\Omega$ is compact, $1 \leq q \leq n$ (\cite{Catlin84b}, \cite{Straube97} when no boundary regularity is assumed). There are natural versions of property($P$) for $(0,q)$-forms, see \cite{FuStraube99}: \eqref{P} is replaced by the requirement that the sum of the smallest $q$ eigenvalues of the complex Hessian of $\lambda_{M}$ should be at least $M$. Note that then $P = P_{1} \Rightarrow P_{2} \Rightarrow \cdots \Rightarrow P_{n}$. This is appropriate, since compactness of $N_{q}$ likewise percolates up the complex: if $N_{q}$ is compact, then so is $N_{q+1}$ (an observation due to McNeal (\cite{McNeal04}), see also the proof of Lemma 2 in \cite{Straube05}). A detailed study of property($P_{1}$) is in \cite{Sibony87b},where various equivalent characterizations are given (with some minimal boundary regularity required for equivalence to the definition we have adopted here, see \cite{FuStraube99}, section 3). In particular, $b\Omega$ satisfies property($P_{1}$) if and only if every continuous function on $b\Omega$ can be approximated uniformly on $b\Omega$ by functions plurisubharmonic in a neighborhood of $b\Omega$. 

In $\overline{\partial}$-problems, the condition that a function be bounded can sometimes be replaced by the condition that the gradient of the function be bounded in the metric induced by the complex Hessian of the function (compare \cite{Berndtsson01} and the references there). In the present context, this was realized by McNeal (\cite{McNeal02}). Say that the boundary of the domain $\Omega$ satisfies condition ($\tilde{P}_{q}$) if, for every $M > 0$, there exists $\lambda_{M} \in C^{2}(U_{M} \cap \Omega)$, where $U_{M}$ is an open neighborhood of $b\Omega$, such that the sum of any $q$ eigenvalues of its complex Hessian is at least $M$, and
\begin{equation}\label{Ptilde}
\left |\sum_{j=1}^{n}\frac{\partial \lambda_{M}}{\partial z_{j}}(z)w_{j}\right |^{2} \leq C\sum_{j,k =1}^{n}\frac{\partial^{2}\lambda_{M}}{\partial z_{j} \partial\overline{z_{k}}}(z)w_{j}\overline{w_{k}}
\end{equation}
for all $w \in \mathbb{C}^{n}$ and $z \in U_{M} \cap \Omega$. That is, $\partial\lambda_{M}$ is bounded uniformly in $M$ in the metric induced by the Hessian of $\lambda_{M}$. McNeal (\cite{McNeal02}) proved the following theorem.
\begin{theorem}\label{ptilde}
Let $\Omega$ be a bounded pseudoconvex domain in $\mathbb{C}^{n}$, $1 \leq q \leq n$. If $\Omega$ satisfies condition($\tilde{P}_{q}$), then $N_{q}$ is compact.
\end{theorem}
\noindent Again note that $\tilde{P}_{1} \Rightarrow \tilde{P}_{2} \Rightarrow \cdots \Rightarrow \tilde{P}_{n}$. Also, there is a weak formulation of \eqref{Ptilde} in terms of currents which is sufficient for Theorem \ref{ptilde}, see \cite{McNeal02} for details.

To prove Theorem \ref{ptilde}, one uses that compactness of $N_{q}$ is equivalent to compactness of the canonical solution operators $\overline{\partial}^{*}N_{q}$ and $\overline{\partial}^{*}N_{q+1}$ (or their adjoints, see e.g. \cite{FuStraube99}, Lemma 1.1). Since $(\tilde{P}_{q}) \Rightarrow (\tilde{P}_{q+1})$, it suffices to establish compactness of $\overline{\partial}^{*}N_{q}$. We sketch an argument (slightly different from McNeal's) to show how the hypotheses of the theorem enter (as well as for use in Remark \ref{Takegoshi} below). We only consider $(0,1)$-forms for simplicity. Also, we assume that $\Omega$ is smooth and that $\lambda_{M}$ is smooth up to the boundary (see \cite{Straube97} and \cite{FuStraube99} for the regularization procedure in the non-smooth case). We may extend $\lambda_{M}$ smoothly to all of $\Omega$ (by shrinking $U_{M}$, where the estimates hold). The starting point is the classical Morrey-Kohn-H\"{o}rmander inequality (\cite{ChenShaw01}, Proposition 4.3.1, and the density Lemma 4.3.2). Taking the weight function to be $\lambda_{M}$ and computing the adjoint of $\overline{\partial}^{*}$ in the weighted metric in  terms of $\overline{\partial}^{*}$ and terms of order zero in $u$ gives for $u \in \Ker(\overline{\partial}) \cap \dom(\overline{\partial}^{*})$
\begin{equation}\label{KoMoHo2}
\int_{\Omega}\sum_{j,k}\frac{\partial^{2}\lambda_{M}}{\partial z_{j} \partial\overline{z_{k}}}u_{j}\overline{u_{k}}e^{-\lambda_{M}} \lesssim
\|\overline{\partial}^{*}(e^{-\lambda_{M}/2}u)\|^{2} +
\|e^{-\lambda_{M} /2}\sum_{j}\frac{\partial\lambda_{M}}{\partial z_{j}}u_{j}\|^{2} \; .
\end{equation}
The constant in \eqref{Ptilde} may be assumed as small as we wish (by scaling $\lambda \rightarrow t\lambda$, see \cite{McNeal02}, p. 199). Therefore, the last term in \eqref{KoMoHo2} can be absorbed into the left hand side, for $z \in U_{M}$. The result is, upon also applying \eqref{P}, using interior elliptic regularity to estimate $\|e^{-\lambda_{M} /2}u\|_{\Omega \setminus \overline{U_{M}}}^{2}\;$ (this term controls the various error terms), and absorbing terms: 
\begin{equation}\label{KoMoHo3}
\|e^{-\lambda_{M} /2}u\|^{2} \lesssim \frac{1}{M}\|\overline{\partial}^{*}(e^{-\lambda_{M} /2}u)\|^{2} + 
C_{M}\|e^{-\lambda_{M} /2}u\|_{-1}^{2} \; , 
\end{equation}
when $u \in \Ker(\overline{\partial}) \cap \dom(\overline{\partial}^{*})$. Denote by $B$ the standard Bergman projection, and by $B_{\varphi}$ the Bergman projection in the weighted space with weight $e^{\varphi}$. Let $v \in \Ker(\overline{\partial})$. Then $v = B(e^{-\lambda_{M} /2}u)$, with $u = B_{-\lambda_{M}/2}(e^{\lambda_{M}/2}v) \in \Ker(\overline{\partial})$. We obtain
\begin{equation}\label{KoMoHo4}
\|v\|^{2} \lesssim \frac{1}{M}\|\overline{\partial}^{*}v\|^{2} + 
C_{M}\|e^{-\lambda_{M}/2}B_{-\lambda_{M}/2}(e^{\lambda_{M}/2}v)\|_{-1}^{2} \; .
\end{equation}
For $M$ fixed, the norm in the last term on the right hand side of \eqref{KoMoHo4} is compact with respect to $\|v\|$, hence with respect to $\|\overline{\partial}^{*}v\|$ (the canonical solution operator to $\overline{\partial}^{*}$ is continuous in $L^{2}$). Since $M$ is arbitrary, this implies that $(\overline{\partial}^{*}N_{1})^{*}$, the canonical solution operator to $\overline{\partial}^{*}$, is compact (see e.g \cite{CDA97}, Lemma 1, \cite{McNeal02}, Lemma 2.1). This concludes the (sketch of) proof of Theorem \ref{ptilde}.

It is easy to see that ${P}_{q}$ implies $\tilde{P}_{q}$, by considering the family $\mu_{M} := e^{\lambda_{M}}$, $M>0$ (see \cite{McNeal02} for details). The exact relationship is however not understood. But there are two classes of domains where property($P_{q}$) is known to be equivalent to compactness of $N_{q}$, and hence condition($\tilde{P}_{q}$) is also equivalent to both compactness of $N_{q}$ and property($P_{q}$). These are the bounded locally convexifiable domains (with no boundary smoothness assumptions assumed beyond what is implied by local convexifiability (Lipschitz), \cite{FuStraube98, FuStraube99}) and the smooth bounded Hartogs domains in $\mathbb{C}^{2}$ (\cite{ChristFu05}). In addition, ($P$) and ($\tilde{P}$) are known to agree for planar domains (\cite{FuStraube02}, Lemma 7). Although condition($\tilde{P}_{1}$) also appears naturally in connection with another important question in several complex variables, namely with that of the existence of a Stein neighborhood basis for the closure of the domain (\cite{Sibony87b}, (4.4) on p. 317), it is not well studied at all. This is in stark contrast to property($P_{q}$), where we have Sibony's theory (\cite{Sibony87b}, see also \cite{FuStraube99} when $q > 1$). It would be very interesting to have an analogous treatment of condition($\tilde{P}_{q}$). Finally, its connection with Stein neighborhoods suggests studying what the implications of compactness in the $\overline{\partial}$-Neumann problem are for the existence of a Stein neighborhood basis for the closure.

Compactness can be viewed as a limiting case of subellipticity. Subellipticity is equivalent to having a bounded plurisubharmonic function, near the boundary, whose Hessian blows up like a power of the reciprocal of the boundary distance (\cite{Catlin87b, Straube97}, compare also \cite{Harrington05, Herbig05}). The only way to show the direction \emph{subellipticity} $\Rightarrow$ \emph{good plurisubharmonic functions} that the author is aware of is via finite type of the boundary. It would be very interesting to have a direct proof of this fact, arguing directly from the subelliptic estimates and bypassing the geometric arguments related to finite type. If there is a reasonable characterization of compactness in terms of (potential theoretic) properties of the boundary (such as $P/\widetilde{P}$ or a slightly weaker property), it should emerge from such a proof when one extracts the features that remain valid in the limiting case.

\begin{remark}\label{Takegoshi} Takegoshi gives a sufficient condition for compactness in \cite{Takegoshi91} which is a precursor to condition($\tilde{P}_{1}$), in the sense that it replaces the uniform boundedness condition on the functions in property($P_{1}$) by a boundedness condition on the gradients. In fact, Takegoshi's condition implies ($\tilde{P}_{1}$), but with \eqref{Ptilde} only for complex tangential $w$'s. But this is enough to prove Theorem \ref{ptilde}, because the forms to which one needs to apply the estimates are in the domain of $\overline{\partial}^{*}$ (see \eqref{KoMoHo2} above), so they are complex tangential (at points of the boundary, but the normal component of a form satisfies a subelliptic estimate, and so terms caused by it are under control). 
\end{remark}
We next want to describe a technique for establishing compactness that does not rely on ($P$)/($\tilde{P}$), introduced recently by the author (\cite{Straube04}). Such a technique is of interest because it is not understood how much stronger (if at all) ($\tilde{P}$) is than compactness. Also, as will be seen, for domains in $\mathbb{C}^{2}$, this technique yields a sufficient condition that modulo a certain (albeit crucial) lower bound on certain radii is also necessary..

Let $\Omega$ be a smooth bounded pseudoconvex domain in $\mathbb{C}^{2}$. If $Z$ is a (real) vector field defined in some open subset of $b\Omega$ (or of $\mathbb{C}^{2}$), we denote by $\mathcal{F}_{Z}^{t}$ the flow genearted by $Z$. In $\mathbb{C}^{2}$, the various notions of finite type coincide (see \cite{D'Angelo93}), so we do not need to specify which notion we mean. Recall that the set of infinite type points in the boundary is compact. Finally, $B(P,r)$ denotes the open ball of radius $r$ centered at $P$. The following theorem is the main result in \cite{Straube04}.
\begin{theorem}\label{04}
Let $\Omega$ be a smooth bounded pseudoconvex domain in $\mathbb{C}^{2}$. Denote by $K$ the set of boundary points of infinite type. Assume there exist constants $C_{1}$, $C_{2}>0$, and $C_{3}$ with $1 \leq C_{3} < 3/2$, and a sequence $\{\varepsilon_{j} > 0\}_{j=1}^{\infty}$ with $\lim_{j \rightarrow \infty}\varepsilon_{j} = 0$ so that the following holds. For every $j \in \mathbb{N}$ and $P \in K$ there is a (real) complex tangential vector field $Z_{P,j}$ of unit length defined in some neighborhood of $P$ in $b\Omega$ with $max|div Z_{P,j}| \leq C_{1}$ such that $\mathcal{F}_{Z_{P,j}}(B(P,C_{2}(\varepsilon)^{C_{3}}) \cap K) \subseteq b\Omega \setminus K$. Then the $\overline{\partial}$-Neumann operator on $\Omega$ is compact.
\end{theorem}
The assumptions in the theorem quantify the notion that at points of $K$, there should exist a (real) complex tangential direction transversal to $K$ in which $b\Omega \setminus K$ (the good set) is `thick' enough. A geometrically very simple special case occurs when $b\Omega \setminus K$ satisfies a complex tangential cone condition (that is, the axis of the cone at a point $P$ lies in a complex tangential direction). In this case, the assumptions in the theorem are satisfied with $C_{3}=1$ (see \cite{Straube04} for details).
\begin{corollary}\label{cone}
Let $\Omega$ be a smooth bounded pseudoconvex domain in $\mathbb{C}^{2}$. Denote by $K$ the set of boundary points of infinite type. Assume that $b\Omega \setminus K$ satisfies a complex tangential cone condition (as an open subset of $b\Omega$). Then the $\overline{\partial}$-Neumann operator is compact.
\end{corollary}
The main idea of the proof of Theorem \ref{04} is very simple. In order to derive a compactness estimate, what needs to be estimated is the $L^{2}$-norm of a form $u$ near $K$. To do this near a point $P$ in $K$, we express $u$ near $P$ in terms of $u$ in a patch which meets the boundary in a relatively compact subset of $b\Omega \setminus K$, plus the integral of the derivative of $u$ in the direction $Z_{P,j}$. The first contribution is easily handled by subelliptic estimates, while the second is dominated by the length of the curve (which is $\varepsilon_{j}$) times the $L^{2}$-norm of $Z_{P,j}u$ on a certain subset of the boundary. But in $\mathbb{C}^{2}$, this norm is dominated by $\|\overline{\partial}u\| + \|\overline{\partial}^{*}u\|$, because $Z_{P,j}$ is complex tangential (so called maximal estimates hold, c.f. \cite{Derridj78}). When one sums up over the various patches (for a fixed $j$), overlap as well as divergence issues arise. These are handled by the uniformity built into the assumptions. 

The conditions in theorem \ref{04} are natural, and, in fact, modulo the size of the lower bound $C_{2}(\varepsilon_{j})^{C_{3}}$ on the radius of the balls, necessary. Indeed, if $N_{1}$ is compact, the boundary contains no analytic discs (since we are in $\mathbb{C}^{2}$, see \cite{FuStraube99} for a proof). This implies, by a result of Catlin (\cite{Catlin78}, Proposition 3.1.12, see also \cite{SahutogluStraube05}, Lemma 3) that for each point $P \in K$ and for every $\varepsilon > 0$, there is a complex tangential vectorfield $Z_{P,\varepsilon}$ (of unit length) near $P$ so that on the integral curve of $Z_{P,\varepsilon}$ through $P$ there is a strictly pseudoconvex point at distance (measured along the curve) less than $\varepsilon$. Then there is a ball $B(P,r)$ which is transported, for $t=\varepsilon$, by the flow generated by $Z_{P,\varepsilon}$, into the points of finite type: it suffices to take $r$ small enough. Since there are smooth bounded pseudoconvex domains in $\mathbb{C}^{2}$ without discs in their boundaries, but whose $\overline{\partial}$-Neumann operator is not compact (\cite{Matheos97}, \cite{FuStraube99}), this discussion also shows that without a lower bound on $r$, the conclusion of the theorem does not hold. The lower bound given in Theorem \ref{04} is probably not optimal. An `optimal' bound (if one exists), in a sense to be made precise, would be of great interest: in light of the above discussion, such a bound essentially amounts to a characterization of compactness in the $\overline{\partial}$-Neumann problem on domains in $\mathbb{C}^{2}$.

Theorem \ref{04} does not hold in dimension $n>2$. Consider a convex domain with a disc in its boundary. When $n>2$, there is an additional complex tangential direction in which to flow, so that the assumptions in Theorem \ref{04} can be satisfied. Yet such domains have noncompact $\overline{\partial}$-Neumann operator (\cite{FuStraube98}). Since the only place where the proof uses that the domain is in $\mathbb{C}^{2}$ is the invocation of maximal estimates, an obvious generalization to $\mathbb{C}^{n}$ is to require the domain to satisfy maximal estimates (equivalently: all the eigenvalues of the Levi form are comparable, see \cite{Derridj78}).There is, however, a more interesting generalization in \cite{Munasinghe06}. It suffices to be able to flow into the set of finite type points along curves whose tangents lie in a complex tangential direction associated with the smallest eigenvalue of the Levi form.

The author does not know examples of domains that satisfy the assumptions in Theorem \ref{04}, but do not satisfy condition($\tilde{P}$). As far as just asserting compactness of the $\overline{\partial}$-Neumann problem on the domains in the theorem is concerned, it does not matter whether or not these domains always satisfy ($\tilde{P}$): we have, in any case, a simple geometric proof of compactness for these domains. However, from the point of view of understanding to what extent ($\tilde{P}$) is necessary for compactness, this question is obviously very important.

\vspace{0.1in}

\textbf{2.2 Obstructions to compactness.} An analytic disc in the boundary constitutes the most blatant violation of  condition($\tilde{P}$) as well as of the condition in Theorem \ref{04}. This is obvious for the condition in Theorem \ref{04} (recall that the setup is in $\mathbb{C}^{2}$). For condition($\tilde{P}$), this can be seen by pulling back the plurisubharmonic functions to the unit disc $D$ in the plane: there do not exist subharmonic functions in $D$ satisfying \eqref{P} and \eqref{Ptilde} for arbitrarily large $M$ (compare Appendix A in \cite{FuStraube99}). Indeed, integration by parts and \eqref{Ptilde} give for $u \in C^{\infty}_{0}(D)$
\begin{multline}\label{unitdisc}
\int_{D}\frac{\partial^{2}\lambda}{\partial z\partial\overline{z}}|u|^{2} = \int_{D}\frac{\partial^{2}\lambda}{\partial z\partial\overline{z}}u\overline{u} \leq
\left |\int_{D}\frac{\partial\lambda}{\partial\overline{z}}\frac{\partial u}{\partial z}\overline{u}\right | +
\left |\int_{D}\frac{\partial\lambda}{\partial\overline{z}}u\frac{\partial \overline{u}}{\partial z}\right | \; \\
\leq \int_{D}\left |\frac{\partial\lambda}{\partial\overline{z}}\right |^{2}|u|^{2} + \int_{D}\left |\frac{\partial u}{\partial z}\right |^{2} \leq
C\int_{D}\frac{\partial^{2}\lambda}{\partial z\partial\overline{z}}|u|^{2} +
\int_{D}\left |\frac{\partial u}{\partial z}\right |^{2} \; ,
\end{multline}
where $C$ is the constant from \eqref{Ptilde}. We have used here that $\|\partial u/ \partial\overline{z}\|^{2} = \|\partial u/ \partial z\|^{2}$. As pointed out earlier, $C$ may be taken as small as we wish. Taking a family with $C = 1/4$ in \eqref{Ptilde} (hence in \eqref{unitdisc}) and combining with \eqref{P} gives 
\begin{equation}\label{false}
\frac{M}{4} \leq  \inf_{u \in C^{\infty}_{0}(D)}\frac{\int_{D}\left |\frac{\partial u}{\partial z}\right |^{2}}{\int_{D}|u|^{2}} \;\, .
\end{equation}
(The infimum on the right hand side of \eqref{false} is (up to a factor $1/4$) the smallest eigenvalue of the Dirichlet realization of $-\Delta$ on $D$.) In the case of property(${P}$), this is of course also an obvious consequence of its characterization (see above) by the approximation property by continuous functions. A disc in the boundary is also known to be an obstruction to hypoellipticity of $\overline{\partial}$ (\cite{Catlin81,DP81}). It is therefore very natural to ask whether such a disc is an obstruction to compactness of the $\overline{\partial}$-Neumann operator. 

An old folklore result, usually attributed to Catlin, says that this is indeed the case for sufficiently regular domains in $\mathbb{C}^{2}$. A proof for the case of Lipschitz boundary may be found in \cite{FuStraube99}. There, a simple example (the unit ball in $\mathbb{C}^{2}$ minus the variety $\{z_{1}=0\}$) is given that shows that some boundary regularity is needed. Whether for domains in $\mathbb{C}^{2}$ there can be other obstructions to compactness was resolved only surprisingly recently. Matheos (\cite{Matheos97}) showed that there are indeed more subtle obstructions: 
\begin{theorem}\label{matheos}
Let $K$ be a compact subset of the complex plane with non-empty fine interior, but empty Euclidean interior. There exists a smooth bounded complete pseudoconvex Hartogs domain in $\mathbb{C}^{2}$ with the following properties: (i) its set of weakly pseudoconvex boundary points projects onto $K$, (ii) it contains no analytic discs in its boundary, (iii) its $\overline{\partial}$-Neumann operator is not compact.
\end{theorem}
\noindent For properties of the fine topology, see e.g. \cite{Helms69}, \cite{FuStraube02}, section 3; in particular, there do exist (many) sets $K$ as in the theorem (an explicit construction of such sets may also be found in section 4 of \cite{ChristFu05}). The version of the theorem given here comes from \cite{FuStraube99}, to where we refer the reader for details (compare also \cite{ChristFu05}). Note in particular that this also means that there are more subtle obstructions to property($P$)/condition($\tilde{P}$) than discs in the boundary. This was known before, see \cite{Sibony87b}.
\begin{remark}\label{hartogs} It is easy to see that on a smooth bounded complete pseudoconvex Hartogs domain in $\mathbb{C}^{2}$, there is no disc in the boundary if and only if the projection of the weakly pseudoconvex boundary points has empty Euclidean interior. It will be seen in subsection 2.3 below that the $\overline{\partial}$-Neumann operator is compact if and only if this set has empty fine interior (as a compact subset of $\mathbb{C}$); see the discussion following Theorem \ref{compandP}. 
\end{remark}
\begin{remark}\label{nodisc}
Whether on a smooth bounded pseudoconvex domain in $\mathbb{C}^{2}$ the absence of discs from the boundary implies global regularity is open.
\end{remark}
It is folklore that the methods that work in $\mathbb{C}^{2}$ can be used in $\mathbb{C}^{n}$ to show that when the $\overline{\partial}$-Neumann operator $N_{1}$ is compact, the boundary cannot contain an $(n-1)$-dimensional complex manifold. However, whether a disc is necessarily an obstruction is open in general, and is arguably the most important problem concerning compactness. \c{S}ahuto\u{g}lu and the author recently showed that when the disc contains a point at which the boundary is strictly pseudoconvex in the directions transverse to the disc (\cite{SahutogluStraube05}), then compactness does fail. This holds more generally for complex submanifolds of the boundary of arbitrary (positive) dimension.
\begin{theorem}\label{disc}
Let $\Omega$ be a smooth bounded pseudoconvex domain in $\mathbb{C}^{n}$, $n \geq 2$. Let $p \in b\Omega$ and assume that the Levi form of $b\Omega$ at $P$ has the eigenvalue zero with multiplicity at most $k$, $1 \leq k \leq n-1$ (i.e. the rank is at least $n-1-k$). If the $\overline{\partial}$-Neumann operator on $(0,1)$-forms is compact, then $b\Omega$ does not contain a $k$-dimensional complex manifold through $P$.
\end{theorem}
It follows immediately from the theorem that if the set $K$ of weakly pseudoconvex boundary points has nonempty relative interior in the boundary, then the $\overline{\partial}$-Neumann operator (on $(0,1)$-forms) is not compact (\cite{SahutogluStraube05}). It suffices to observe that near a relative interior point of $K$ where the Levi form attains its maximal rank (among points of the relative interior of $K$), say $m$, the rank has to be constant, so that the boundary is foliated there by complex manifolds of dimension $n-1-m$. Note that in general, $K$ is considerably bigger than the set of Levi \emph{flat} points.

The proof of Theorem \ref{disc} results from the following ideas. Compactness is a local property (see \cite{FuStraube99}, Lemma 1.2). Therefore, it suffices to argue locally. Assume the boundary contains a complex manifold, say $M$. There is a holomorphic change of coordinates near $P$ so that in the new coordinates $M$ is affine, and the real normal to the boundary of $\Omega$ is constant on $M$. This is always possible (\cite{SahutogluStraube05}, Lemma 1). Next, consider a section $\Omega_{1}$ of $\Omega$ through $P$, perpendicular to $M$. If $\Omega_{1}$ has a subdomian $\Omega_{2}$, whose boundary shares $P$ with $b\Omega$ and such that (i) the restriction operator from the Bergman space of $\Omega_{1}$ to the Bergman space of $\Omega_{2}$ is not compact, and (ii) the product $\Omega_{2} \times M$ is contained in $\Omega$ (near $P$), then the arguments from \cite{FuStraube98} (which in turn are based on ideas from \cite{Catlin81, DP81}) carry over to produce a contradiction to the existence of a compact solution operator to $\overline{\partial}$ (which would be a consequence of compactness of $N_{1}$). In the situation of Theorem \ref{disc}, any smooth subdomain $\Omega_{2}$ will do, because $\Omega_{1}$ is strictly pseudoconvex at $P$ (\cite{SahutogluStraube05}, Lemma 2). If we take for $\Omega_{2}$ a ball with small radius and tangent to $b\Omega_{1}$ at $P$, (ii) also holds (because the real normal to $b\Omega$ is constant along $M$).

Experience indicates that a flatter boundary should be even more favorable to noncompactness of the $\overline{\partial}$-Neumann operator. In other words, the extra assumption that the boundary is strictly pseudoconvex in the directions transverse to $M$ should not be needed. However, the present methods do not seem to yield this. An interesting recent contribution to this circle of ideas, involving the Kobayashi metric of the domain, is in \cite{Kim04}.

The above proof of Theorem \ref{disc} also raises a question of independent interest. Namely given a domain $\Omega$ and a subdomain $\Omega_{1}$, when is the restriction operator from the Bergman space of $\Omega$ to that of $\Omega_{1}$ compact? Of course, this happens when $\Omega_{1}$ is relatively compact in $\Omega$, so the case of interest is that where the domains share a boundary point. As mentioned above, it is also known to happen when $\Omega$ is smooth near a strictly pseudoconvex boundary point $P$ and $b\Omega_{1}$ shares $P$ with $b\Omega$ and is smooth there (\cite{SahutogluStraube05}, Lemma 2). In addition, this restriction is known to be compact when $\Omega$ is convex, $P=0 \in b\Omega$, and $\Omega_{1} = r\Omega$, for $r<1$ (\cite{FuStraube98}). The general situation is not understood.

\vspace{0.1in}

\textbf{2.3 Hartogs domains in $\mathbb{C}^{2}$ and semi-classical analysis of Schr\"{o}dinger operators.} It is well known that $\overline{\partial}$ and related operators on Hartogs domains can be studied by means of weighted operators on the base domain. In the sequel, the base domain $U$ will be a planar domain (i.e. the Hartogs domain is in $\mathbb{C}^{2}$). The resulting weighted problems lead to Schr\"{o}dinger operators on $U$, see for example \cite{Berndtsson96} and the references there.

Let $U$ be a bounded domain in $\mathbb{C}$, $\phi(z) \in C^{2}(\overline{U})$. Denote by $S_{\phi}$ the Schr\"{o}dinger operator with magnetic potential $A = -(\partial\phi / \partial y)dx + (\partial\phi / \partial x)dy$, magnetic field $dA = \Delta\phi (dx\wedge dy)$, and electric potential $V = \Delta\phi$. That is, $S_{\phi}$ is given by (the Dirichlet realization of)
\begin{equation}\label{schr}
S_{\phi} = -\left [(\partial /\partial x + i\partial\phi / \partial y)^{2} +
(\partial / \partial y - i\partial\phi / \partial x)^{2} \right ] + 
\Delta \phi \; . 
\end{equation}
Denote by $S_{\phi}^{0}$ the corresponding nonmagnetic Schr\"{o}dinger operator, given by (the Dirichlet realization of)
\begin{equation}\label{schr2}
S_{\phi}^{0} = -\Delta + \Delta\phi \; .
\end{equation}
For (very) brief introductions to Schr\"{o}dinger operators, we refer the reader to \cite{FuStraube02}, section 2 or \cite{ChristFu05}, section 2.3. For a detailed treatment in the context of semi-classical analysis relevant here, see \cite{Helffer88}

Let $\Omega$ be a bounded complete pseudoconvex Hartogs domain in $\mathbb{C}^{2}$ given by $\Omega = \{(z,w) \in \mathbb{C}^{2}: z \in U, |w| < e^{-\phi(z)}\}$, where $U$ is a domain in $\mathbb{C}$. Note that pseudoconvexity forces $\phi$ to be plurisubharmonic. Also, smoothness of $\Omega$ means that $\phi$ is smooth on $U$, but not on $\overline{U}$, but the notions needed here are still well defined, compare \cite{FuStraube02}, \cite{ChristFu05}. The rotation invariance in the $w$ variable brings a discrete Fourier variable into play, and so what one actually has when analyzing the $\overline{\partial}$ and related problems on Hartogs domains are sequences of Schr\"{o}dinger operators of the form $\{S_{n\phi}\}_{n=1}^{\infty}$ and $\{S_{n\phi}^{0}\}_{n=1}^{\infty}$, respectively (see \cite{FuStraube02} for details). Compactness of the $\overline{\partial}$-Neumann operator on $\Omega$ is closely linked to the behavior of the sequence of lowest eigenvalues $\{\lambda_{n\phi}\}_{n=1}^{\infty}$ (the ground state energies) of the magnetic Schr\"{o}dinger operators, while property($P$) of $b\Omega$ is similarly linked to the behavior of the sequence $\{\lambda_{n\phi}^{0}\}_{n=1}^{\infty}$ of lowest eigenvalues of their nonmagnetic counterparts. The former idea originates with \cite{Matheos97}, the latter with \cite{FuStraube02}. The precise relationships are given in the following theorem (\cite{FuStraube02}).
\begin{theorem}\label{equiv}
Let $\Omega = \{(z,w) \in \mathbb{C}^{2}: |w| < e^{-\phi(z)}, z \in U\}$ be a smooth bounded complete pseudoconvex Hartogs domain. Suppose that $b\Omega$ is strictly pseudoconvex on $b\Omega \cap \{w=0\}$. Then
\vspace{0.05in}

\noindent (1) $b\Omega$ satisfies property ($P$) if and only if $\lambda_{n \phi}^{0} \rightarrow \infty$ as $n \rightarrow \infty$. 
\vspace{0.04in}

\noindent (2) The $\overline{\partial}$-Neumann operator on $\Omega$ is compact if and only if $\lambda_{n\phi} \rightarrow \infty$.
\end{theorem}
\noindent We remark that for the domains in Theorem \ref{equiv}, property ($P$) and property ($\tilde{P}$) are equivalent (\cite{FuStraube02}, Lemma 6). It was already noted in \cite{FuStraube02} that for some of the implications, regularity of the boundary is not needed. For a version of Theorem \ref{equiv} that assumes very little regularity of $\phi$, see \cite{ChristFu05}.

It is the limit in part $(2)$ of Theorem \ref{equiv} that gives rise to the terminology used in the title of this subsection. Note that $S_{n\phi} = -n^{2}[\left ((1/n)(\partial /\partial x) + i(\partial \phi /\partial y)\right )^{2} + \left ((1/n)(\partial /\partial y) -i(\partial\phi /\partial x)\right )^{2}] + n\Delta\phi $. Understanding the behavior of the ground state energy as $n$ tends to infinity is thus analogous to understanding (modulo the factor $n^{2}$) what happens when `Planck's constant' $h = 1/n$ tends to zero. This situation is typically referred to as semi-classical analysis in the mathematical physics literature. Mathematical physics also has its own version of $(1) \Rightarrow (2)$: Simon's diamagnetic inequality asserts that $\lambda_{n\phi}^{0} \leq \lambda_{n\phi}$ (\cite{Simon76}, see also \cite{Kato72}). Reverse relationships, when there is some kind of domination of the magnetic eigenvalues by the nonmagnetic ones, obviously of interest in our context, are known in the  physics literature as paramagnetism. For more thorough discussions of these topics, we refer again to \cite{Helffer88}, \cite{FuStraube02}, and \cite{ChristFu05}, and their references.

This point of view has allowed to clarify the relationship between property($P$) and compactness of the $\overline{\partial}$-Neumann operator on the (special) class of Hartogs domains in $\mathbb{C}^{2}$. Namely, Christ and Fu recently established the paramagnetic property required for the implication $(2) \Rightarrow (1)$ in Theorem \ref{equiv}.
\begin{theorem}\label{christfu}
Let $\phi$ be subharmonic on the domain $U \subseteq \mathbb{C}$, and suppose that $\Delta\phi$ is H\"{o}lder continuous of some positive order. If    $\;\sup_{n}\lambda_{n\phi}^{0} < \infty$ then $\liminf_{n \rightarrow \infty}\lambda_{n\phi}$  $< \infty$.
\end{theorem}
Since the sequence $\{\lambda_{n\phi}^{0}\}_{n=1}^{\infty}$ is increasing (this is obvious from \eqref{schr2}; $\phi$ is subharmonic), we immediately get the corollary that on the domains from Theorem \ref{equiv}, property($P$) and compactness of the $\overline{\partial}$-Neumann operator are equivalent. In fact, combining this with some additional work, Christ and Fu (\cite{ChristFu05}) were able to handle general (not necessarily complete) Hartogs domains, thus establishing the following equivalence.
\begin{theorem}\label{compandP}
Let $\Omega \subseteq \mathbb{C}^{2}$ be a smooth bounded pseudoconvex Hartogs domain. The $\overline{\partial}$-Neumann operator on $\Omega$ is compact if and only if $b\Omega$ satisfies property($P$).
\end{theorem}
While we have not made an effort to state Theorem \ref{compandP} with optimal boundary smoothness assumptions (see \cite{ChristFu05}), we point out that in view of the example mentioned before the statement of Theorem \ref{matheos}, some boundary regularity is needed for the equivalence in Theorem \ref{compandP} to hold. If the Hartogs domain is complete, then the two properties in Theorem \ref{compandP} are also equivalent to the set $K$ (as defined above) having empty fine interior, by work of Sibony. Namely, $b\Omega$ satisfies property($P$) if and only if $K$ does (as a subset of $\mathbb{C}$, \cite{Sibony87b}, p. 310). In turn, $K$ satisfies property($P$) if and only if it has empty fine interior (\cite{Sibony87b}, Proposition 1.11).

\section{Global regularity} 

The $\overline{\partial}$-Neumann operator $N_{q}$ is said to be globally regular if it maps $C^{\infty}{(0,q)}(\overline{\Omega})$ (necessarily continuously) into itself. It is said to be exactly regular if it maps $W^{s}_{(0,q)}(\Omega)$ into itself for $s \geq 0$. Exact regularity implies of course global regularity. So far, in all instances where one can prove global regularity, one actually proves exact regularity. On the worm domains, failure of global regularity (\cite{Christ96}) is proved via failure of exact regularity (\cite{Barrett92}): for most $s$, exact a priori estimates hold in $W^{s}(\Omega)$, and global regularity would then give exact regularity. It is consistent with what is known that such a priori estimates might hold on all domains (they are also known to hold on the nonpseudoconvex counterexample domains from \cite{Barrett84}, see \cite{BoasStraube92}).  

For a survey of results up to about ten years ago, we refer he reader to \cite{BoasStraube99}. In this section, we first discuss regularity on domains whose boundary contains an open patch foliated by complex hypersurfaces (\cite{StraubeSucheston03, ForstnericLaurent05}. In subsection 3.2, we describe a unified approach to global regularity (\cite{StraubeSucheston02, Straube05}). 

We recall the following important $1$-form on the boundary of a domain. Let $\Omega$ be a smooth bounded pseudoconvex domain. Denote by $\eta$ a purely imaginary nowhere vanishing $1$-form on $b\Omega$ that annihilates the complex tangent space and its conjugate. Let $T$ denote the purely imaginary vector field on $b\Omega$ orthogonal to the complex tangent space and its conjugate and such that $\eta(T) \equiv 1$. The real $1$-form $\alpha$ is defined by $\alpha = -\mathcal{L}_{T}\eta$, the Lie derivative of $\eta$ in the direction of $T$ (compare \cite{D'Angelo87, D'Angelo93}). The form arises naturally in the computation of (normal components of) commutators of vector fields; indeed, if $\eta = \partial\rho - \overline{\partial}\rho$, and $\overline{X}$ is a local section of $T^{0,1}(b\Omega)$, then $\alpha(\overline{X}) = 2\partial\rho([L_{n}, \overline{X}])$ ($L_{n}$ is the complex normal). The cohomology class on complex submanifolds of the boundary mentioned in the introduction in connection with \cite {BoasStraube93} is the class of $\alpha$.

\vspace{0.1in}
\textbf{3.1 A foliation in the boundary.}  Background on foliation theory and notions used here can be found in \cite{CC2000} and \cite{Tondeur97}, as well as in \cite{StraubeSucheston03}, \cite{ForstnericLaurent05}, and their references. Assume now there is a codimension one foliation in the boundary, say the relative interior of the set $K$ of weakly pseudoconvex points is foliated by complex manifolds of dimension $n-1$. Note that such a foliation is always transversely orientable (by the vector field $T$ defined on all of $b\Omega$). In order to run the machinery from \cite{BoasStraube91, BoasStraube93}, one needs a function $h$ smooth in a relative neighborhood of $K$, satisfying 
\begin{equation}\label{leaf}
dh|_{L} = \alpha|_{L} \; , \;\mathrm{for\ all\ leaves}\ L \; .
\end{equation}
For details, see \cite{StraubeSucheston03}. Of course, this requires that the restriction of $\alpha$ to a leaf is closed. This does indeed hold: $d\alpha |_{\mathcal{N}_{P}} = 0$ always, where $\mathcal{N}_{P}$ is the null space of the Levi form at $P$, see the lemma in section 2 of \cite{BoasStraube93}. Thus solving \eqref{leaf} is always possible locally.  Globally, topological constraints arise. Also, the boundary behavior of $h$ on $K$ needs to be controlled.

It turns out that these issues are very much related to ones studied in foliation theory. Note that the foliation can be defined by $\eta$: the tangent planes to the leaves are given by the null space of $\eta$. Then the Frobenius condition reads $d\eta\wedge\eta = 0$. Hence $d\eta = \beta\wedge\eta$ for some $1$-form $\beta$. $\alpha$ is such a form, that is
\begin{equation}\label{d-eta}
d\eta = \alpha\wedge\eta \; \;\mathrm{on} \;K \;  
\end{equation}
(\cite{Tondeur97}, Proposition 2.2). With \eqref{d-eta}, solving \eqref{leaf} is easily tied to an important property in foliation theory (\cite{StraubeSucheston03}).
\begin{lemma}\label{closed}
\eqref{leaf} can be solved (say on the relative interior of $K$) if and only if the Levi foliation of $K$ can be defined \emph{globally} by a \emph{closed} $1$-form.
\end{lemma}
\noindent Indeed, if $\omega$ is a $1$-form defining the foliation, then $\omega = e^{-h}\eta$, and 
\begin{equation}\label{d-omega1}
d\omega = d(e^{-h}\eta) = e^{-h}(-dh\wedge\eta + d\eta) = e^{-h}(-dh + \alpha)\wedge\eta \; .
\end{equation}
Therefore, 
\begin{equation}\label{d-omega2}
d\omega = 0 \Leftrightarrow (-dh + \alpha)\wedge\eta = 0 \Leftrightarrow -dh|_{L} + \alpha|_{L} = 0 \; . 
\end{equation}

In addition to closedness of $\alpha|_{L}$, solvability of \eqref{leaf} also requires that the De Rham cohomology class of $\alpha|_{L}$ vanishes. This again fits nicely into the foliation framework: this cohomology class coincides with the infinitesimal holonomy of $L$ (\cite{StraubeSucheston03}, Remark 2, \cite{CC2000}, Example 2.3.15).

We first present a result in $\mathbb{C}^{2}$ from \cite{StraubeSucheston03}. The relative boundary of $K$ in $b\Omega$, say $\Gamma$, is assumed smooth, and so is a smooth compact orientable surface embedded in $\mathbb{C}^{2}$. Recall that a  complex tangency at a point of $\Gamma$ is called generic if it is either elliptic or hyperbolic (see \cite{StraubeSucheston03} for more information). Note that at a hyperbolic point there are locally two leaves that meet. 
\begin{theorem}\label{foliation1}
Let $\Omega$ be a smooth bounded pseudoconvex domain in $\mathbb{C}^{2}$. Suppose that the set $K$ of infinite type points of $b\Omega$ is smoothly bounded (in $b\Omega$) and that its boundary $\Gamma$ is connected and has only isolated generic complex tangencies. Assume that the two leaves meeting at a hyperbolic point are distinct globally and that they have no other hyperbolic points in their closure (in $K$). If each leaf of the Levi foliation is closed (in the relative interior of $K$) and has trivial infinitesimal holonomy, then the $\overline{\partial}$-Neumann operator on $\Omega$ is continuous on $W_{(0,1)}^{s}(\Omega)$ for $s \geq 0$.
\end{theorem}
If one assumes that the Levi foliation of $K$ is part of a foliation of a bigger smooth Levi flat hypersurface $M$, with $M \cap \overline{\Omega} = K$, then boundary behavior of the leaves is easier to control, and one needs no conditions on the boundary of $K$. This results in a geometrically appealing sufficient condition. The following theorem is from \cite{ForstnericLaurent05}. A codimension one foliation is called simple if through every point there exists a local transversal (a line) that meets each leaf at most once.
\begin{theorem}\label{foliation2}
Let $\Omega \subseteq \mathbb{C}^{n}$ be a smooth bounded pseudoconvex domain such that the set $K$ of all boundary points of infinite D'Angelo-type is the closure of its relative interior in $b\Omega$. Assume $K$ is contained in a smooth Levi-flat (open) hypersurface $M \subset \mathbb{C}^{n}$, whose Levi foliation satisfies one (hence both) of the following equivalent conditions: (i) the leaves of the restriction of the foliation to a neighborhood of $K$ are topologically closed (ii) the foliation is simple in a neighborhood of $K$. Then, the $\overline{\partial}$-Neumann operators $N_{q}$ on $\Omega$ are continuous in $W^{s}_{(0,q)}(\Omega)$ for $s \geq 0$, $1 \leq q \leq n$.
\end{theorem}
Very roughly speaking, in both Theorems \ref{foliation1} and \ref{foliation2}, one would like to obtain a closed form that defines the foliation by pulling back from the leaf space (which is one dimensional) a form roughly like $dx$. One then has to deal with the non-Hausdorff nature of this space. In Theorem \ref{foliation1}, one also has to control the boundary behavior.

It is interesting to note that the main concern in \cite{ForstnericLaurent05} is not the $\overline{\partial}$-Neumann problem, but rather holomorphic convexity properties of compact subsets of $M$. The authors use asymptotically pluriharmonic defining functions for $M$ (near a compact subset) for constructing Stein neighborhoods, and whether such defining functions exist leads precisely to the question whether the foliation, near the compact set (globally), can be defined by a closed $1$-form. In view of Lemma \ref{closed}, this is related to the equivalence between `pluriharmonic defining functions' and `exactness of $\alpha$ in \cite{StraubeSucheston02} (see Theorem \ref{equivalent} below). In a local context, compare also \cite{BF88}.

The two equivalent conditions in Lemma \ref{closed} are equivalent to a third one, given in terms of the flow generated by $T$ (\cite{StraubeSucheston03}, Proposition 2). This is at least potentially of interest because the condition is in terms of $T$, (rather than the Levi foliation), which is well defined on the boundary of any smooth domain. Furthermore, this leads to a homological necessary and sufficient condition for the existence of a function $h \in C^{\infty}(K)$ as above (\cite{StraubeSucheston03}, Theorem 3), in terms of foliation currents (\cite{Sullivan76, CC2000}) associated to $T$. 

\vspace{0.1in}
\textbf{3.2 Sufficient conditions for global regularity.} In \cite{BoasStraube93}, Remark 3 in section 4, the authors point out that the families of vector fields with good approximate commutator conditions with $\overline{\partial}$, required in their approach to global regularity (\cite{BoasStraube91, BoasStraube93}) can exist in situations where the domain does not admit (even a local) plurisubharmonic defining function. On the other hand, they had noted in \cite{BoasStraube91} that it suffices to have the commutator conditions with components of $\overline{\partial}$ in directions that lie in the null space of the Levi form. For this, positivity of the Hessian of a defining function at a boundary point is needed only on the span of the null space of the Levi form and the complex normal. The situation was cleared up in \cite{StraubeSucheston02}: the authors showed that the vector fields and the plurisubharmonic defining functions approaches can be reformulated naturally and then become equivalent.

Let $\Omega$ be a smooth bounded pseudoconvex domain. Say that $\Omega$ admits a family of essentially pluriharmonic defining functions if there exists a family $\{\rho_{\varepsilon}\}_{\varepsilon > 0}$ of defining functions with gradients bounded and bounded away from zero on $b\Omega$ uniformly in $\varepsilon$, such that the complex Hessian of $\rho_{\varepsilon}$ is $O(\varepsilon)$ on the span, over $\mathbb{C}$, of $\mathcal{N}_{P}$ and $L_{n}(P)$, for all $P \in b\Omega$. We emphasize that this notion is indeed a generalization of the notion of a pluri\emph{sub}harmonic defining function (see \cite{StraubeSucheston02}). We say that the form $\alpha$ (see above) is approximately exact on the null space of the Levi form if there exists a family $\{h_{\varepsilon}\}_{\varepsilon > 0}$ of functions smooth in neighborhoods $U_{\varepsilon}$ of the set $K$ of boundary points of infinite D'Angelo type, bounded uniformly in $\varepsilon$, such that $dh|_{\mathcal{N}_{P}} = \alpha|_{\mathcal{N}_{P}} + O(\varepsilon)$ for all $P \in K$. A family of conjugate normals which are approximately holomorphic in weakly pseudoconvex directions is defined similarly; see \cite{StraubeSucheston02}, where Sucheston and the author established the following equivalence (compare also the remarks in section 5 of \cite{Straube05} concerning $(iii)$).
\begin{theorem}\label{equivalent}
Let $\Omega$ be a smooth bounded pseudoconvex domain in $\mathbb{C}^{n}$. The following are equivalent:

(i) $\Omega$ admits a family of essentially pluriharmonic defining functions.

(ii) $\Omega$ admits a family of conjugate normals which are approximately holomorphic in weakly pseudoconvex directions.

(iii) $\Omega$ admits a family of vector fields as in \cite{BoasStraube93}.

(iv) The form $\alpha$ is approximately exact on the null space of the Levi form.
\end{theorem}
The equivalence to condition (ii) is of interest because the existence of such a family leads, under favorable circumstances ($K$ is uniformly H-convex), to the existence of transverse vector fields \emph{holomorphic in a neighborhood of $K$}, and these lead to Stein neighborhood bases for $\overline{\Omega}$ and to Mergelyan type approximation (\cite{FN77}, \cite{StraubeSucheston02}, \cite{ForstnericLaurent05}).
 
Theorem \ref{equivalent} shows that the approaches to global regularity in the $\overline{\partial}$-Neumann problem through plurisubharmonic defining functions and through good vector fields are really equivalent. Left unanswered was the question of how to unify this approach with that via compactness. This is achieved in the following theorem from \cite{Straube05}.
\begin{theorem}\label{unify}
Let $\Omega$ be a smooth bounded pseudoconvex domain in $\mathbb{C}^{n}$, $\rho$ a defining function for $\Omega$. Let $1 \leq q \leq n$. Assume that there is a constant $C$ such that for all $\varepsilon > 0$ there exists a defining function $\rho_{\varepsilon}$ for $\Omega$ and a constant $C_{\varepsilon}$ with
\begin{equation}\label{uniform}
1/C < |\nabla\rho_{\varepsilon}| < C \ \ \mathrm{on}\  b\Omega \; ,
\end{equation}
and
\begin{equation}\label{sufficient}
\left \|\sum^{\prime}_{|K|=q-1} \left (\sum_{j,k =1}^{n}\frac{\partial^{2}\rho_{\varepsilon}}{\partial z_{j}\partial\overline{z_{k}}}\frac{\partial\rho}{\partial\overline{z_{j}}}
\overline{u_{kK}} \right ) d\overline{z_{K}} \right \|^{2} \leq
\varepsilon \left (\|\overline{\partial}u\|^{2} + \|\overline{\partial}^{*}u\|^{2} \right ) + C_{\varepsilon}\|u\|_{-1}^{2}
\end{equation}
for all $u \in C^{\infty}_{(0,q)}(\overline{\Omega}) \cap \dom(\overline{\partial}^{*})$. Then 
\begin{equation}\label{estimate}
\|N_{q}u\|_{s} \leq C_{s}\|u\|_{s} \; ,
\end{equation}
for $s \geq 0$ and all $u \in W^{s}_{(0,q)}(\Omega)$.
\end{theorem}
Notice that the assumptions in Theorem \ref{unify} are for $q$-forms, $q$ fixed. It is not hard to see that when they are satisfied at level $q$, then they are satisfied at level $q+1$ (\cite{Straube05}, Lemma 2). It would be interesting to know whether global regularity similarly moves up to higher form levels (recall from section 2 that subellipticity and compactness do).

The simplest situation occurs when there is one defining function, say $\rho$, that works for all $\varepsilon$. This covers the case when $N_{q}$ is compact: the left hand side is in this case bounded by $\|u\|^{2}$ independently of $\varepsilon$, and compactness says precisely that $\|u\|^{2}$ can be bounded in the manner required by the right hand side of \eqref{sufficient} (this is the right hand side of a compactness estimate).

When $\Omega$ admits a defining function $\rho$ that is plurisubharmonic at the boundary, $\rho_{\varepsilon} = \rho$ for all $\varepsilon$ also works. Assume $q=1$ for the moment. Applying the Cauchy-Schwarz inequality to the left hand side of \eqref{sufficient} at boundary points gives that this left hand side is dominated by $\|\sum(\partial^{2}\rho / \partial z_{j} \partial\overline{z_{k}})u_{j}\overline{u_{k}}\|^{2}$ plus a term of order $\rho$ plus a compactly supported term. The latter two are benign for \eqref{sufficient}. Estimating the former in the way required in \eqref{sufficient} is immediate from the fact that the Hessian of a defining function acts as a subelliptic multiplier of order $1/2$ on $1$-forms (\cite{D'Angelo93}, section 6.4.2). When $q > 1$, one can reformulate \eqref{sufficient} so that the left hand side of the inequality involves a pairing between $q$-forms (\cite{Straube05}, Lemma 1), and the above argument works under the weaker assumption that the sum of any $q$ eigenvalues of the Hessian of $\rho$ is nonnegative. In view of the equivalence results in \cite{BoasStraube90}, this recovers, in the pseudoconvex case, a recent result of Herbig-McNeal (\cite{Herbig-McNeal05}), where the authors prove Sobolev estimates for the Bergman projection on $j$-forms, $q-1 \leq j \leq n$, under this weaker assumption.

More generally, the sufficient conditions for global regularity from Theorem \ref{equivalent} imply those in Theorem \ref{unify}:
\begin{proposition}\label{unify2}
Let $\Omega$ be a smooth bounded pseudoconvex domain in $\mathbb{C}^{n}$. Assume $\Omega$ satisfies one (hence all) of the equivalent conditions in Theorem \ref{equivalent}. Then the assumptions in Theorem \ref{unify} are satisfied for $q = 1, 2,  \cdots , n$.
\end{proposition}
We indicate what is involved, details are in \cite{Straube05}, Proposition 1. Assume $(i)$ in Theorem \ref{equivalent}. It suffices to consider the case $q=1$. Fix $\varepsilon$. Then $L_{\rho_{\varepsilon}}(z)(w) = O(\varepsilon)|w|^{2}$ when $(z,w)$ is in a neighborhood $U_{\varepsilon}$ of the compact subset $\{(z,w): w \in \mathcal{N}_{z}\}$ of the unit sphere bundle in $T^{1,0}(b\Omega)$. Here, $L_{g}$ denotes the complex Hessian of a function $g$. There is a constant $C_{\varepsilon}$ such that $|w|^{2} \leq C_{\varepsilon}L_{\rho_{\varepsilon}}(z)(w)$ when $(z,w) \notin U_{\varepsilon}$. This implies the estimate, when $z \in b\Omega$, $w \in T^{1,0}(b\Omega)$:
\begin{equation}\label{split}
\left |\sum_{j,k=1}^{n}\frac{\partial^{2}\rho_{\varepsilon}}{\partial z_{j}\partial\overline{z_{k}}}(z)\frac{\partial\rho}{\partial\overline{z_{j}}}(z)\overline{w_{k}} \right |^{2} \leq
C\varepsilon|w|^{2} +
\widetilde{C_{\varepsilon}}\sum_{j,k=1}^{n}\frac{\partial^{2}\rho_{\varepsilon}}{\partial z_{j}\partial\overline{z_{k}}}(z)w_{j}\overline{w_{k}} 
\end{equation}
(since both terms on the right are nonnegative). By continuity and homogeneity, \eqref{split} holds near (depending on $\varepsilon$) the boundary. To verify \eqref{sufficient} for $u \in C^{\infty}_{(0,1)}(\overline{\Omega})$, it suffices to apply \eqref{split} to $u$ pointwise, near the boundary (the normal component of $u$ is zero only on the boundary, but it satisfies a subelliptic estimate, so is under control). In view of the discussion preceeding the statement of Proposition \ref{unify2}, integration over $\Omega$ now gives \eqref{sufficient}.

It should not be surprising that condition \eqref{sufficient} has a potential theoretic flavor: global regularity probably is not determined by geometric conditions alone (unlike the much stronger property of subellipticity). However, it is not hard to extract a geometric sufficient condition from \eqref{sufficient}, compare \cite{Straube05}, section 2. What one arrives at is precisely condition $(i)$ in Theorem \ref{equivalent}. In other words, \emph{the vector field approach constitutes} what might be called \emph{the geometric content of Theorem \ref{unify}}.

It is noteworthy that whether or not a family of defining functions satisfies \eqref{uniform} and \eqref{sufficient} is determined entirely by the interplay of the gradients with the boundary. That is, if a family $\{\rho_{\varepsilon}\}_{\varepsilon > 0}$ satisfies \eqref{uniform} and \eqref{sufficient}, and $\{\widetilde{\rho_{\varepsilon}}\}_{\varepsilon > 0}$ is another family such that $\nabla(\widetilde{\rho_{\varepsilon}}) = \nabla(\rho_{\varepsilon})$ for all $\varepsilon$ and all $z \in b\Omega$, then $\{\widetilde{\rho_{\varepsilon}}\}_{\varepsilon > 0}$ also satisfies \eqref{uniform} and \eqref{sufficient} (possibly after rescaling). For details, see \cite{Straube05}, Remark 2.

At the level of a priori estimates, a proof of Theorem \ref{unify} follows from a small modification of the ideas in \cite{BoasStraube91, BoasStraube93}.  We briefly indicate what changes, keeping the general setup from \cite{BoasStraube91}. This will show how \eqref{sufficient} enters into the estimates. Set $X_{\varepsilon} = e^{-h_{\varepsilon}}\sum(\partial\rho / \partial \overline{z_{j}})(\partial / \partial z_{j})$, where $h_{\varepsilon}$ is defined by $\rho_{\varepsilon} = e^{h_{\varepsilon}}\rho$, and $\rho$ is a defining function with normalized gradient (all of this is near $b\Omega$, away from $b\Omega$, any smooth continuation will do). For the Bergman projection $P$, the key quantity to be estimated is
\begin{equation}\label{mainterm}
\left \|\varphi(X_{\varepsilon} - \overline{X_{\varepsilon}})Pf \right\|^{2} 
\lesssim
\left (N_{1}\overline{\partial}f, \varphi^{2}(X_{\varepsilon} - \overline{X_{\varepsilon}})[\overline{\partial}, X_{\varepsilon} - \overline{X_{\varepsilon}}]Pf \right ) + \mathrm{o.k.} \; ,
\end{equation}
where `o.k.' stands for terms that are under control or can be absorbed, and  $\varphi$ is a smooth cutoff function supported near the boundary (see \cite{BoasStraube91}, p. 83--84). One needs to control the normal component of the commutator in \eqref{mainterm}. In contrast to \cite{BoasStraube91}, we do not have a pointwise estimate on this normal component. But there is some slack built into the argument in \cite{BoasStraube91}, in that there the contribution from the commutator of $ X_{\varepsilon}-\overline{X_{\varepsilon}}$ with each component of $\overline{\partial}$ is estimated separately. If one takes this into account, computing the commutator gives the main term (after integrating $X_{\varepsilon}-\overline{X_{\varepsilon}}$ back to the left)
\begin{equation}\label{mainterm2}
\left (\sum_{j,k=1}^{n}\frac{\partial^{2}\rho_{\varepsilon}}{\partial z_{j}\partial\overline{z_{k}}}\left ((X_{\varepsilon} - \overline{X_{\varepsilon}})N_{1}\overline{\partial}f\right )_{j}\frac{\partial\rho}{\partial z_{k}}\;, \;X_{\varepsilon}Pf \right )
\end{equation}
(as opposed to estimating $\sum_{k=1}^{n}\cdots$ for each $j$, $j=1, \cdots, n$). The term in the left hand side of this inner product is now (the conjugate of) one to which \eqref{sufficient} can be applied. (As usual, we let $X_{\varepsilon} - \overline{X_{\varepsilon}}$ act in special boundary charts so that it preserves $\dom(\overline{\partial}^{*})$.) Note that $X_{\varepsilon} = e^{h_{\varepsilon}}L_{n}$, and that $e^{h_{\varepsilon}}$ is bounded independently of $\varepsilon$. Consequently (by \eqref{sufficient}), the square of the $L^{2}$-norm of this term is dominated by
\begin{multline}\label{final}
\varepsilon \left (\|\overline{\partial}(L_{n}-\overline{L_{n}})N_{1}\overline{\partial}f\|^{2} +
\|\overline{\partial}^{*}(L_{n}-\overline{L_{n}})N_{1}\overline{\partial}f\|^{2} \right ) \\ + C_{\varepsilon}\|(L_{n}-\overline{L_{n}})N_{1}\overline{\partial}f\|_{-1}^{2} \\
\lesssim \varepsilon \left (\|N_{1}\overline{\partial}f\|_{1}^{2} + 
\|\overline{\partial}^{*}N_{1}\overline{\partial}f\|_{1}^{2} \right ) + 
C_{\varepsilon}\|f\|^{2} \; .
\end{multline}
From here on, the argument proceeds as in \cite{BoasStraube91}; in particular, in the setup of the downward induction on $q$ there, $N_{1}\overline{\partial}$ is `as good as' $P$. Absorbing terms, one arrives at the required a priori estimate (compare p. 84--85 in \cite{BoasStraube91}).

The situation changes rather markedly with regard to genuine estimates. In \cite{BoasStraube91}, the authors simply observe that the estimates can be carried out uniformly on suitable approximating strictly pseudoconvex subdomains, by using the same family of vector fields. By contrast, the assumptions in Theorem \ref{unify} do not seem strong enough to be inherited by these approximating subdomains. Therefore, one has to employ some other regularization procedure, such as elliptic regularization. This makes the argument considerably more involved, and the author derives in \cite{Straube05} certain needed new estimates for the regularized operators. This is also in contrast to \cite{ChenShaw01}, where the results of \cite{BoasStraube91} are proved working directly with the $\overline{\partial}$-Neumann operator and using elliptic regularization. There too it is the strength of the pointwise estimates on the size of the normal component of the commutators that makes elliptic regularization routine (once the derivation of the a priori estimates is in place).

Note that to get estimates at a fixed Sobolev level $k$, it suffices to have \eqref{sufficient} in Theorem \ref{unify} for some $\varepsilon = \varepsilon(k)$. Kohn (\cite{Kohn99}) has proved estimates where the level in the Sobolev scale up to which estimates hold is tied to the Diederich-Forn\ae ss exponent (\cite{DF77b}) of the domain. The discussion above of Proposition \ref{unify2}, combined with \cite{Straube01}, where the plurisubharmonicity of $-\log(-\rho)$ is exploited, suggests that it should be possible to obtain results of this type by the methods in \cite{Straube05}.
\begin{remark}\label{multiplier}
Consider the operator $A_{\rho}$ from  $\dom(\overline{\partial}) \cap \dom(\overline{\partial}^{*})$, provided with the graph norm, to $L^{2}(\Omega)$, given by
\begin{equation}\label{sesqui}
A_{\rho}(u) = \sum_{j,k=1}^{n}\frac{\partial^{2}\rho}{\partial z_{j}\partial\overline{z_{k}}}\frac{\partial\rho}{\partial\overline{z_{j}}}\overline{u_{k}} \; , \; u \in \dom(\overline{\partial}) \cap \dom(\overline{\partial}^{*})  \; . 
\end{equation}
Then \eqref{sufficient} holds with $\rho_{\varepsilon} = \rho$ for all $\varepsilon$ precisely when $A_{\rho}$ is compact (see e.g. \cite{CDA97}, Lemma 1, \cite{McNeal02}, Lemma 2.1). The form of $A_{\rho}$ suggests that one study sesquilinear forms that produce compact operators via \eqref{sesqui}. It is possible that there is a theory of `compactness multipliers'. We mention that compactness of $A_{\rho}$ for a suitable defining function $\rho$ is considerably weaker than compactness of $N_{1}$. It holds on all convex domains (since they admit a plurisubharmonic defining function), yet $N_{1}$ is compact (if and) only if the boundary of the domain contains no analytic disc (\cite{FuStraube98}).
\end{remark}

\section{References}

\renewcommand{\refname}{}    
\vspace*{-36pt}              

\providecommand{\bysame}{\leavevmode\hbox to3em{\hrulefill}\thinspace}

\end{document}